# Mathematical Pressure Volume Models
# of the Cerebrospinal Fluid

by

S. Sivaloganathan and G. Tenti
Department of Applied Mathematics
University of Waterloo
Waterloo, Ontario, Canada N2L 3G1

and

J.M. Drake
Division of Neurosurgery
Hospital for Sick Children
University of Toronto
Toronto, Ontario, Canada M5G 1X8

ABSTRACT

Numerous mathematical models have emerged in the medical literature over the past two decades attempting to characterize the pressure and volume dynamics of the central nervous system compartment. These models have been used to study the behavior of this compartment under such pathological clinical conditions as hydrocephalus, head injury and brain edema. The number of different approaches has led to considerable confusion regarding the validity, accuracy or appropriateness of the various models.

In this paper we review the mathematical basis for these models in a simplified fashion, leaving the mathematical details to appendices. We show that most previous models are in fact particular cases of a single basic differential equation describing the evolution in time of the cerebrospinal fluid pressure (CFS). Central to this approach is the hypothesis that the rate of change of CSF volume with respect to pressure is a measure of the compliance of the brain tissue which as a consequence leads to particular models depending on the form of the compliance function.

All such models in fact give essentially no information on the behavior of the brain itself. More recent models (solved numerically using the Finite Element Method) have begun to address this issue but have difficulties due to the lack of information about the mechanical properties of the brain. Suggestions are made on how development of models which account for these mechanical properties might be developed.

**Key Words:**    Pressure-Volume relation. CSF dynamics. Hydrocephalus.

1. INTRODUCTION

In this paper we present a comprehensive and critical review of the mathematical models which have been proposed in the medical literature over the last twenty years in order to gain insight into the hydrodynamical behavior of the brain cerebrospinal fluid (CSF).

It is not our intention to review here the extensive literature on the physiology of the CSF, for which we refer the reader to well-known classical monographs [1] and more recent compilations [2,3]. Instead, we shall focus our attention on a detailed presentation of the mathematical models of CSF dynamics proposed thus far for two main reasons: Firstly, we believe with Hakim [4] that, before attempting to interpret every aspect of intracranial physiopathology on the basis of biological phenomena alone, it is very important to examine such problems in the light of classical concepts of physics. And, secondly, mathematical modeling of such phenomena has the virtue of allowing one to make precise predictions from the provisional hypotheses of the model, which in turn is very helpful in the design of critical experiments and in the interpretation of their results.

The application of non-trivial mathematical concepts and techniques to the study of CSF hydrodynamics is a relatively recent development. Of course, there have been theoretical studies before; but they were either concerned with special problems [5], or they failed to demonstrate the relevance of the mathematical model to the biological situation [6]. In 1972, however, Guinane [7] published a mathematical analysis of his experimental results, based on a simple differential equation (DE) for the CSF pressure, which can be clearly recognized as the first non-trivial mathematical model in this field.

One of the central assumptions of Guinane's model is that the CSF pressure $P(t)$ at any time $t$ is linearly related to the CSF volume $V(t)$ at that same time. Although for small changes of volume and pressure this linear relation is supported by experiment, there is clear evidence that for larger deviations from equilibrium the so-called "pressure-volume relation" is more nearly exponential. This led to the next generation of mathematical models a few years later, when Marmarou **et al.** [8] proposed a nonlinear DE for the pressure, in which the exponential form of the pressure-volume relation was incorporated, and used it to study three cases of clinical interest, namely bolus injection, volume removal, and infusion. This model quickly became one of the best known and stimulated further work along similar lines, such as Drake's extension to include the presence of shunts and the concomitant negative pressure-volume relationship [9].

The next class of models appeared much more recently in the work of Nagashima **et al.** [17], and represents a substantial increase in the level of mathematical sophistication. In contrast to the phenomenological models, Nagashima **et al.** aim at reaching the much more ambitious goal of deriving the effective stress distribution in the brain along with the interstitial pressure change and the ventricular configuration. Their starting point is Hakim's assumption [13] that the brain parenchyma may be schematized as a microscopic sponge of viscoelastic material, and their mathematical model is a first approximation of Hakim's concept of "open cell sponge" in an attempt to describe mathematically the mechanical interaction between the brain parenchyma and the CSF.

The approach suggested by Nagashima **et al.** [17] relies heavily on the theory of consolidation introduced by Biot [5] half a century ago for the treatment of certain types of flow through porous media, and well known in soil mechanics. It is a continuum theory and leads to a set of partial differential equations which are usually solved numerically by standard techniques such as the Finite Element Method. As usual in problems of continuum mechanics, one of the crucial aspects of the theory is the discussion of which type of boundary conditions are appropriate for the problem at hand. Unfortunately, Nagashima **et al.** [17] do not give the reader any details of the reasoning which led to their choice of boundary conditions, with the result that it is very difficult to make an assessment of their work at the present time. Undoubtedly, it represents the right path to explore in the future; but it would be highly desirable to elucidate the underpinnings of consolidation theory as applied to the brain before embarking on numerical solutions of the mathematical model.

## 2. PHENOMENOLOGICAL MODELS

On the basis of many studies carried out over several decades, we consider the space-averaged CSF pressure $P$ as an essential variable. Furthermore, to a first approximation the fluid may be considered incompressible (like water) and confined to a single compartment of volume $V$. The next assumption comes from the experimental finding that the CSF is produced at a certain rate $I_f$ and absorbed at a certain rate $I_a$, in locations which are immaterial for the present argument. Since both of these rates may vary in time, it follows that both the volume of the CSF space $V(t)$, and the fluid pressure $P(t)$, are functions of time as well.

In order to derive an equation for the evolution of the pressure we start by applying the principle of mass conservation, according to which the rate of change of the volume must match the difference between the rates of formation and absorption, viz.

$$\begin{aligned} \frac{dV}{dt} &= \text{(rate of CSF formation)} - \text{(rate of CSF absorption)} \\ &= I_f(t) - I_a(t) \ . \end{aligned} \qquad (1)$$

Next we invoke the results of infusion studies to postulate that the absorption rate is given by

$$I_a(t) = \frac{1}{R_a}\left[P(t) - P_d\right], \qquad (2)$$

where the constant $R_a$ is a fluid resistance, and the threshold pressure $P_d$, above which the absorption mechanism operates, is usually taken to be the sagittal sinus venous pressure. The rate of formation, instead, does not depend on pressure, and is assumed to be known once the specific experimental design has been chosen.

The final step of the derivation consists of expressing the left-hand side of (??) in terms of the pressure. Because of the assumed incompressibility of the CSF, it is clear that an increase of the fluid's pressure will produce an increase of the volume of the CSF space by a compression of the brain parenchyma. Consequently, it seems reasonable to postulate that the volume $V$ depends on time through the intermediacy of the pressure $P$; mathematically, this means that $V(t) = V(P(t))$, and it follows that

$$\begin{aligned} \frac{dV}{dt} &= \frac{dV}{dP}\frac{dP}{dt} \ , \\ &= C(P)\frac{dP}{dt} \ , \end{aligned} \qquad (3)$$

where

$$C(P) = \frac{dV}{dP} \ , \qquad (4)$$

which is referred to as the brain compliance. Substitution of (??) and (??) into (??) and simple re-arrangement of the terms gives the desired equation for the pressure, namely

$$C(P)\frac{dP}{dt} + \frac{P(t)}{R_a} = I_f(t) + \frac{P_d}{R_a} \ , \qquad (5)$$

which is a first-order DE for the function $P$, with initial condition

$$P(0) = P_0. \qquad (6)$$

It is interesting to note that (??) represents not just one mathematical model but a whole family of models, each of which is specified by the choice of a particular form for the compliance $C(P)$.

We shall review the various choices which have been made in the literature in the next sections, and since the compliance $C(P)$ is not derived from fundamental principles but rather inferred from experimental results, we shall refer to all such models as phenomenological models.

**Linear Phenomenological Models.**

The simplest choice for the compliance is that it be appreciably independent of pressure, so that we may assume that $C(P) = C_0 =$ constant, as done by Guinane [12] on the grounds that it appears to be experimentally verified over a limited range of pressures. Next we must specify the quantity $I_f(t)$, which represents the combined effects of the CSF formation rate in the brain and the amount of artificial CSF injected or withdrawn at the initial stage.

Although the detailed mechanism of CSF formation is quite complicated [8] it is not too important for our primitive model. All we require at this stage is the information that, under normal physiological conditions, a state of equilibrium (or steady state) exists; thus it appears reasonable to assume a constant rate of formation $I_f^{(e)}$ – i.e. independent of pressure – and hence to express the quantity $I_f(t)$ in the form

$$I_f(t) = I_f^{(e)} + S(t) \quad , \quad \left( I_f^{(e)} = \text{constant} \right) \tag{7}$$

where $S(t)$ is an external source specified in concrete cases by the experimentalist's procedure (e.g. bolus injection). When these assumptions are used, we obtain our first class of models (of which Guinane's is a member) given by:

$$C_0 \frac{dP}{dt} + \frac{P(t)}{R_a} = I_f^{(e)} + \frac{P_d}{R_a} + S(t) , \tag{8}$$

which must be solved subject to the initial condition given by (??). The complete mathematical solution to this problem is presented in Appendix A. However, for the reader who is interested only in the results, we offer a synopsis in Table 1. The first model is the one originally due to Guinane [12], with the source function set equal to zero. As the solution shows, once the CSF pressure has been somehow disturbed away from the equilibrium value

$$P_{\text{eq}} = R_a I_f^{(e)} + P_d \tag{9}$$

it relaxes back to it with a frequency $\nu_a = \dfrac{1}{C_0 R_a}$ (i.e. a relaxation time $\tau = C_0, R_a$). Thus, this simplest of models requires four input parameters: $C_0, R_a, I_f^{(e)}$, and $P_d$, which must be obtained from experimental measurements. A graphical representation of this solution is given in Fig. 1,

where - as in all other numerical simulations shown in this paper which are illustrative only - we have used the typical values of the parameters listed in Table 2.

Model #2 of Table 1 is obtained from the basic equation (??) by assuming that the external source $S(t)$ acts instantaneously. For example, one may produce a pressure pulse by injecting or withdrawing from the ventricular space a certain amount of CSF. The procedure takes, of course, a finite time interval; however, it is relatively easy to perform it in such a way that this time interval may be considered as infinitely short in comparison to the typical relaxation time of the system.

Here $S_0$ is a constant which is taken to be positive if CSF is added (bolus injection) or negative if CSF is subtracted (volume withdrawal), while $t_0$ represents the instant at which the pulse is produced.

The solution shown in the fifth column consists of two parts, the first of which is identical with Model #1 while the second one contains the Heaviside, or unit step, function

$$H(t - t_0) = \begin{cases} 1 & \text{if } t \geq t_0, \\ 0 & \text{otherwise}. \end{cases} \tag{10}$$

merely reflecting the fact that the pressure pulse is felt at the instant $t_0 > 0$ at which the source or sink is switched on. The size of the pressure pulse externally produced is given by the new parameters $\overline{P} = \dfrac{S_0}{C_0}$ which can be positive or negative, depending on the sign of $S_0$. In the first case (illustrated in Fig. 2) starting with an initial pressure $P_0 > P_{\text{eq}}$ arbitrarily chosen, we see the CSF pressure decays towards the equilibrium value as before. However, when the bolus injection at $t = t_0$ produces the pressure pulse $\overline{P}$ we observe a sudden jump to the value $P(t_0) + \overline{P}$ in the CSF pressure, which slowly falls back thereon to the final equilibrium value as $t \to \infty$. In the second case, (illustrated in Fig. 3) the sudden withdrawal of CSF at $t = t_0$ results in a sudden decrease, which is followed thereafter by a slow recovery to the final equilibrium value $P_{\text{eq}}$.

Model # 3 in Table 1 represents a third situation of clinical interest, namely that of continuous infusion at a rate $S_i$ started at a selected time instant $t_i$. Then we can represent the source function in the form shown in the third column, where the Heaviside function $H$ guarantees that there is no infusion until the chosen time $t_i$. The solution shown in the fifth column consists again of two parts, the first of which is exactly the same as the previous models. The second part, on the other hand, contains, contains the new parameters $P_i = R_i S_i$ (where $R_i$ stands for the fluid resistance of the infusion apparatus) and $\nu_i = \dfrac{1}{C_0 R_i}$ which is the corresponding frequency. This part represents a very different behavior, reflecting the different experimental conditions, for it shows that after a long time from the beginning of the infusion procedure a new steady state value of the pressure is reached given by

$$P_{ss} = P_{\text{eq}} + \frac{\nu_i}{\nu_a} P_i . \tag{11}$$

Physically, and within the assumptions of this class of models, this means that before infusion starts the pressure relaxes from the initial value $P_0$ in the same manner as Model # 1, that is to say with a time constant $\tau_a = C_0 R_a$ determined by the value of the brain compliance and the fluid resistance to absorption. Immediately after the start of the infusion process, in contrast, the volume of the CSF is rapidly increased and so the CSF pressure experiences a corresponding increase; this in turn produces an increase in the absorption rate, according to (??), and so this model predicts that the brain adjusts by reaching the new equilibrium value $P_{ss}$ which will persist as long as infusion continues. The behavior of the CSF pressure just described is vividly illustrated by Fig. 4.

The final model in Table 1 (Model #4) includes the effect of a shunt on the CSF dynamics. It is well known that shunts are a reasonably effective treatment of hydrocephalus, and that they work by diverting CSF from the ventricular system or lumbar theca into another anatomical compartment where fluid can be absorbed. From the standpoint of this paper, the presence of a shunt in the CSF circulation system may be easily taken into account by modifying the phenomenological model represented by (??) to read

$$\frac{dV}{dt} = (\text{rate of CSF formation}) - (\text{rate of CSF absorption}) \\ - (\text{rate of shunting}). \tag{12}$$

For purposes of illustration, we consider the example of a standard one-way valve whose flow rate is described (in the simplest possible terms) by the equation

$$I_s(t) = \frac{1}{R_s} [P(t) - P_{op}] \, H(t - t_s) . \tag{13}$$

Here $R_s$ is the valve resistance to flow, and the parameter $P_{op}$ is a characteristic of the shunt and is normally specified by the manufacturer. The constant $t_s$ appearing in the Heaviside function may be assumed to represent the time necessary to reach a situation where $P(t_s) > P_{op}$, so that the shunt starts operating.

¿From a strictly mathematical point of view, very little modification of our general linear model (??) is required; in fact, all that is necessary is the replacement of the external source function $S$ by the function $I_s$ defined in (??). After the introduction of the shunt frequency $\nu_s = \dfrac{1}{C_0 R_s}$, the differential equation takes on the form shown in the fourth column of Table 1, and its solution for times longer than $t_s$ is displayed in the fifth column. For the sake of brevity a sketch of the graph of this function is not shown here, since it is very similar to that shown in Fig. 1 – except, of course, for the fact that the equilibrium value is different, being given by

$$P^* = \frac{\nu_a P_{\text{eq}} - \nu_s P_{op}}{\nu_a - \nu_s} , \qquad (14)$$

rather than $P_{\text{eq}}$. More complicated shunts require a separate analysis, and are beyond the scope of this review.

**Nonlinear Phenomenological Models.**

The various models analyzed above were all predicated on the basic assumption of a constant compliance. As mentioned already in the introduction, however, there exists strong experimental evidence that this is not true over a large interval of pressures, so that a better class of models is desirable. A step in this direction was taken in 1978 by Marmarou **et al.** [16] who postulated a compliance function of the form

$$C(P) = \frac{1}{kP} , \qquad k = \text{constant} \qquad (15)$$

This assumption implies that the pressure-volume relation is exponential in nature, as can easily be seen from (??) which now becomes

$$\frac{dV}{dP} = \frac{1}{kP} , \qquad (16)$$

whose solution is just

$$P(V) = P' e^{k(V - V')} , \qquad (17)$$

where $V'$ is some reference volume corresponding to the reference pressure $P'$. Denoting the volume differences by $\Delta V$, dividing through by the reference pressure $P'$ and taking the natural logarithm of both sides of (??) gives

$$\ln\left(\frac{P}{P'}\right) = k \, \Delta V , \qquad (18)$$

or else $\frac{1}{k} = \Delta V / \ln\left(\frac{P}{P'}\right)$, which goes under the name of the pressure-volume index (PVI). Thus a measurement of the PVI allows one to find a numerical value for the constant $k$.

This model of brain compliance is listed as Model # 2 in Table 3 and is a favorite among clinicians and researchers alike, as already mentioned, not only because it fits the data reasonably well but also because of its extreme simplicity. In fact, it is clear form (??) that in a semi-logarithmic plot (i.e., volume versus the logarithm of the pressure) the PVI is just the slope of the resulting straight line.

A slight modification of the Marmarou **et al.** model was suggested by Sklar and his collaborators [20,21] and is listed as Model #3 in Table 3. It used two parameters instead of one and, as a result, it has the more realistic feature of predicting a finite compliance when the pressure vanishes, instead of the infinite value predicted by Model # 2. For $P > 0$, however, the graph of $C(P)$ is still given by a branch of a hyperbola, and the PVI is given by

$$\frac{1}{k_1} = \frac{\Delta V}{\ln\left(\frac{P + \frac{k_2}{k_1}}{P' + \frac{k_2}{k_1}}\right)}, \tag{19}$$

which is basically the same as the one defined by Marmarou **et al.**

Finally, we have listed as Model #4 the one suggested by Lim **et al.** [15] in which the brain compliance is assumed to decrease exponentially with increasing CSF pressure. These authors found that with an appropriate choice of the parameters $a$ and $b$ they could fit their experimental results on normal and hydrocephalic dogs quite well, and in particular that the value of the parameter $a$ appears to be strongly dependent on whether the animals were normal or hydrocephalic, while the value of $b$ appeared to be the same. We would like to point out, however, that this exponential model does not appreciably add to the hyperbolic model proposed by Marmarou **et al.**, in the sense that they both fit the available data well, as Sklar and his collaborators have repeatedly pointed out [20,21]. This is to be expected since noting that the parameter $b$ is very small [16] and that the range of values taken by the pressure is limited, we may write this model in the following form

$$C(P) = \frac{a}{e^{bP}} \simeq \frac{a}{1 + bP}, \tag{20}$$

using the Maclaurin expansion of the exponential truncated after the second term. The resulting approximation is just the hyperbolic model, as anticipated, and we will not consider Model #4 any further.

The rest of this section is devoted to a thorough analysis of the DE resulting from the assumption that the brain compliance follows the hyperbolic law postulated by model #2. (The use of Model #3 leads to an entirely similar set of results.) The reason for this detailed discussion is that although this model was devised by Marmarou **et al.** [16] over fifteen years ago, and is by now considered a classic, we must point out a few flaws in these authors' analysis.

When assumption (??) for the compliance is used in our general phenomenological Model (??), and the formation rate is written as in (??), the result is a nonlinear, first-order DE, viz.

$$\frac{dP}{dt} - \alpha(P_{\text{eq}} + R_a S(t)) P + \alpha P^2 = 0 \tag{21}$$

where $\alpha = \dfrac{k}{R_a}$ and the remaining symbols have the same meaning as before. This is a Riccati equation about which much is known (see [14]).

It should be noted that, as was the case with its linear counterpart (??), equation (??) defines a whole class of nonlinear models, each of which is determined by a particular choice of the source function $S$. However, rather than solving anew each case, it is possible and convenient to derive the general solution of this DE for an arbitrary source function, as done in Appendix B, and then to discuss the various cases. Hence we now move on to a discussion of their physical meaning as it follows from the mathematical analysis presented in the appendix.

*Case 1.    No external source.*

The situation here is similar to the linear case, that is to say the initial pressure $P_0$ is taken to be above the equilibrium value $P_{\text{eq}}$ and the subsequent temporal behavior is observed. The solution of (??) is given in Appendix B (B6) and is repeated here for the reader's convenience:

$$P(t) = \frac{P_{\text{eq}}}{1 + \left(\frac{P_{\text{eq}}}{P_0} - 1\right) e^{-\alpha P_{\text{eq}} t}} \ . \tag{22}$$

A comparison of this with the corresponding result of the linear theory (Model #1 of Table 1) shows that the two models differ only in one respect, namely that the nonlinear model predicts a somewhat faster relaxation to equilibrium, as illustrated in Fig.5 for typical values of the parameters. In particular, it is noteworthy that different forms for the brain compliance only affect the transient behavior of the pressure, but not its final equilibrium value. This is the result of the fact that only the relaxation time is changed, from $\tau_a = R_a C_0$ for the linear model to $\tau_\alpha = \alpha P_{\text{eq}}$ for the nonlinear one.

*Case 2. Bolus injection (or withdrawal).*

We can represent the external stimulus by a combination of Heaviside functions. To illustrate, let $t_b$ be the instant at which the bolus injection is started, and suppose that it lasts a time $\Delta t_b$ during which an amount $S_b$ of fluid is injected per unit time into the CSF space; then the source function $S(t)$ may be written as

$$S(t) = S_b \left[ H(t - t_b) - H(t - t_b - \Delta t_b) \right] , \tag{23}$$

i.e. $S(t) = S_b$ for $t_b \leq t \leq t_b + \Delta t_b$ and vanishes otherwise. When $S_b > 0$ the graph of this function is illustrated in Fig. 6, while for $S_b < 0$ (bolus withdrawal) it is given by the mirror image with respect to the time axis.

The representation of $S$ by (??) makes the coefficient of $P(t)$ in the differential equation (??) a discontinuous function; but this represents no problem from the mathematical point of view. As shown in Appendix B, we simply solve the differential equation over the three separate intervals $0 \leq t < t_b$, $t_b \leq t < t_b + \Delta t_b$, and $t_b + \Delta t_b \leq t < \infty$, and then match the solution at the points $t_b$ and $t_b + \Delta t_b$. Albeit a little complicated, those formulas are the ones needed when one is interested in describing in detail what happens to the CSF pressure before, during, and after the bolus injection is performed.

Suppose, however, that we were not interested in the behavior of the pressure before and during the bolus injection, but only in its relaxation to equilibrium after the disturbance has been produced. Then the mathematical problem could be significantly simplified by proceeding as follows. Starting from the equilibrium pressure $P_{\text{eq}}$ we increase the CSF volume by the amount $\Delta V = S_b \Delta t_b$ through a bolus injection; then, according to the assumed exponential relation (??), the new CSF pressure is given by

$$P(\Delta V) = P_{\text{eq}} e^{k S_b \Delta t_b} . \tag{24}$$

As long as the duration of the procedure $\Delta t_b$ is much shorter than the relaxation time of $P(t)$, we may assume that this value of the pressure is reached instantaneously, after which the external source is no longer present. In other words, we solve (??) with $S(t) = 0$ and take the equation above as the initial value written more simply as

$$P_0 = P(0) = P_{\text{eq}} \, e^{\nu'_b \Delta t_b} , \tag{25}$$

where the notation introduced in Appendix B has been used. Using the procedure explained therein produces the following result:

$$P(t) = \frac{P_{\text{eq}}}{1 + \left(e^{-\nu_b' \Delta t_b} - 1\right) e^{-(\nu_b - \nu_b')t}} \qquad (26)$$

Aside from the notation, this formula is exactly the same as formula I-D in Fig. 2 of Marmarou **et al.** [16]; however, their interpretation of it in Fig. 2 I-F is misleading, for it seems to describe the transient behavior of $P$ before and during the bolus injection as well as afterwards. As explained above, (??) is only valid **after** the pressure is raised to the initial value $P_0$ given by (??), with the consequence that its graph resembles that of our Fig. 5, and not that of Fig. 2 of [16].

Similar remarks apply to the case of a volume removal, which differs only in having $S_b < 0$; it follows that $\nu_b'$ is also negative and hence the initial pressure $P_0$ given by (??) is smaller than in the previous case. As a result, the recovery of the CSF pressure to the value $P_{\text{eq}}$, which it had before the volume removal, appears typically as in Fig. 7, and not as in Fig. 2 II-F of [16] which is only correct in the part which they labelled $b \leftrightarrow c$. Presumably, these authors were attempting to include the case in which the removal of a volume $\Delta V$ lowers the CSF pressure so much that the absorption mechanism ceases to operate. But if that were the case then the assumption made in (??) – which is the same as equation (2) of [16] – for the absorption rate $I_a(t)$ would not be valid, and hence the basic differential equation (??) (or equation (7) of [16]) would cease to be valid and more mathematical care must be exercised, i.e. we could simply take care of that possibility by replacing (??) with

$$I_a(t) = \frac{1}{R_a}[P(t) - P_d]\, H(t - t_a), \qquad (27)$$

where $t_a$ is the instant at which CSF absorption starts to occur, and the presence of the Heaviside function shuts off absorption until that instant. The corresponding differential equation then reads

$$\frac{1}{kP}\frac{dP}{dt} = I_f(t) - \frac{1}{R_a}[P(t) - P_d]\, H(t - t_a). \qquad (28)$$

For the case of interest here, that is a bolus withdrawal of short duration $\Delta t_b$ which brings the initial pressure down to the initial value

$$P_0 = P_{\text{eq}} e^{-\nu_b' \Delta t_b} < P_{\text{eq}}, \qquad (29)$$

(??) reduces to the simple form

$$\frac{dP}{dt} = k I_f^{(e)} P, \quad (0 \leq t < t_a) \qquad (30)$$

whose solution subject to (??) is given by

$$P(t) = P_{\text{eq}} e^{-\nu'_b \Delta t_b + k I_f^{(e)} t}, \quad (0 \leq t < t_a) \tag{31}$$

which coincides with the result of Marmarou **et al.** [16] (given by formula (1) in Fig. 2, case II-D).

Finally, we consider the case of continuous infusion at a constant rate for which the source function is given in Model #3 of Table 1. Assuming that the absorption mechanism is in continuous operation, the nonlinear equation (??) now reads

$$\frac{dP}{dt} - \alpha \left[ P_{\text{eq}} + R_a S_i H(t - t_i) \right] P + \alpha P^2 = 0, \tag{32}$$

which we shall solve subject to the initial condition $P(0) = P_{\text{eq}}$ in order to compare our result with that of Marmarou **et al**. Before infusion starts, i.e. over the time interval $0 \leq t \leq t_i$, the initial-value problem reduces to

$$\frac{dP}{dt} - \alpha P_{\text{eq}} P + \alpha P^2 = 0,$$

$$P(0) = P_{eq}, \tag{33}$$

whose only solution is $P(t) = \text{constant} = P_{eq}$, as the reader can verify immediately. On the other hand, after infusion is started and continued thereafter we have $H(t - t_i) = 1$, and so the initial-value problem becomes

$$\frac{dP}{dt} - \alpha (P_{eq} + R_a S_i) P + \alpha P^2 = 0,$$

$$P(0) = P_{eq} \tag{34}$$

whose solution is simply given by

$$P(t) = \frac{P_{eq} + R_a S_i}{1 + \frac{R_a S_i}{P_{eq}} e^{-\alpha(P_{eq} + R_a S_i)t}}, \quad (t \geq t_i) \tag{35}$$

which coincides with formula III-D in Fig. 2 of Marmarou **et al.** [16]. On taking the limit as $t \to \infty$ we see from (??) that the pressure reaches a new steady-state value given by

$$P_{ss} = P_{eq} + R_a S_i. \tag{36}$$

Moreover, by introducing the "resistance to infusion" $R_i$ and writing

$$R_a S_i = \frac{R_a}{R_i}(R_i S_i) = \frac{\nu_i}{\nu_a} P_i \ , \tag{37}$$

as we did for the linear model, we see that the $P_{ss}$, predicted by the nonlinear model in (**??**) exactly coincides with that predicted by the linear model in (**??**). Thus the different functional form for the brain compliance does not affect the final equilibrium value reached by the pressure.

An analysis completely analogous to the one just given can be done for the case in which shunts are taken into account. However, we shall spare the reader such lengthy calculations and simply state that the behavior found is very similar to the linear case already discussed and that, in fact, the equilibrium pressure is exactly the same as (**??**) of the linear case (for the same type of shunt). All the details may be obtained from the solution in Table 1 (Model #4).

3. DISCUSSION

Several phenomenological models describing the transient behavior of the CSF pressure have been developed over the past two decades, and these have helped in the understanding of the pathological changes in head injury, hydrocephalus, and brain edema. We have attempted to present these mathematical models from a unified point of view, by deriving a basic differential equation (**??**), from which the various models can be constructed by specifying the functional form of the so-called brain compliance. Thus, by assuming a constant compliance one obtains a class of linear models of which Guinane's [12] is a typical example, while the more realistic assumption of a pressure dependent compliance leads to nonlinear models such as those put forward by Marmarou **et al.** [16], or by Drake **et al.** [9].

An interesting feature of all these models, quite independent of the mathematical complications, is their prediction of a common steady-state value for the CSF pressure. As we have shown, drastic changes in the functional form of the compliance lead to different relaxation times, and hence to relatively small differences in the transient behavior of the pressure; but they do not affect in any way its final equilibrium value. ¿From a purely mathematical point of view this result looks peculiar, since, in general, nonlinear differential equations show a much richer and varied behavior than their linear counterparts - particularly so, insofar as the existence of equilibria and their stability are concerned.

It is therefore worthwhile taking a critical look at the fundamental assumptions of the phenomenological models. As explained earlier (Cf. (**??**) and (**??**)), the basic differential equation (**??**) rests on the hypothesis that the rate of change of the CSF volume with respect to the pressure $\frac{dV}{dP}$ is a measure of the compliance of the brain tissue, including blood vessels. The reasoning behind it is that , the total cranium volume being a constant (at least in the adult), any increase in the CSF volume must be produced by a compression of the brain tissue. In other words, the function $C(P)$ defined by (**??**) is considered a measure of the volume storage capacity determined by the elastic properties of the system, and as such is affected by the compressibility of the vascular volume as well as the subpial brain tissue. Plausible as this may sound, the fact remains that one should first of all show that there is indeed a correlation between the function $C(P)$ and the actual elastic (or viscoelastic) properties of the subpial brain tissue. Not only has this not been done mathematically, but there exists recent experimental evidence showing that no such correlation exists [25].

It is important to appreciate, of course, that these negative results do not invalidate Hakim's basic picture of the brain as a porous medium, for there is no doubt that the extra volume made available to the CSF is obtained at the expense of the solid matrix. What they do indicate, however, is that the quantity $\frac{dV}{dP}$ is not an appropriate parameter for obtaining information about the brain's elastic properties. On the other hand, the pressure-volume response is easily measured, and it would be clearly desirable to tie this parameter to one or more specific properties of the brain tissue. Several attempts have been made in this direction, notably by Benabid and his collaborators [3] who defined the compliance as the rate of change of the cerebral blood volume with respect to intracranial pressure, that is to say as the vascular bed elasticity. Unfortunately, in recognition of the difficulty of formulating a satisfactory theory of the phenomena involved, these authors recommend that attention be focussed on the use of experimental data to extract information on the parameter $\frac{K}{S}$, which is the ratio of the CSF pathways conductance to vascular elasticity. The result is that their phenomenological models for the bolus injection procedure, as well, as for the case of constant-rate injection, are simply members of the family of linear models reviewed earlier. Thus the possible connections between the CSF pressure and the material properties of the brain still remain an open problem, deserving of a great deal of attention in future investigations in this area.

A second important direction for future research is the development of mathematical models built on Hakim's picture of the brain as a porous medium. As mentioned in the Introduction, this program has been recently started by Nagashima **el al.** [17] who integrated numerically the equations of the theory of consolidation from Soil Mechanics equations (**??**)-(**??**). The latter was developed to explain the settlement of a soil under load, which is caused by a gradual adaptation of the soil to the load variation. A very successful theory of this one-dimensional system was given long ago by Terzaghi [22,23], and its limitations removed by Biot later on [4,5]. As a result, the modern theory of consolidation is based on the following system of partial differential equations

$$\epsilon = \nabla \cdot \vec{u}, \tag{38}$$

$$\mu \nabla^2 \vec{u} + (\lambda + \mu) \nabla \epsilon = \nabla \sigma, \tag{39}$$

$$\frac{K}{\rho g} \nabla^2 \sigma = \frac{\partial \epsilon}{\partial t} + n\beta' \frac{\partial \sigma}{\partial t}, \tag{40}$$

where $\nabla$ and $\nabla^2$ are the gradient and the Laplacian operators, respectively. The vector $\vec{u}$ represents the incremental displacement vector of the solid skeleton, $\epsilon$ the incremental volume strain, and $\sigma$ the incremental fluid pressure. The coefficients $\lambda$ and $\mu$ are the Lamé coefficients of the theory of elasticity, $g$ is the acceleration due to gravity, $K$ the so-called permeability and $\rho$ the fluid density. Finally $n$ stands for the porosity of the medium and $\beta'$ for the fluid's compressibility.

Equations 38-40 form the basic mathematical model on which the numerical computations of Nagashima **et al.** [17] are based. Valuable as these calculations are, it would seem highly desirable to submit this model to detailed analytical studies in order to gain insight into its region of applicability as well as its limitations. For example, the hypothesis under which Biot derived these equations are repeated by these authors, but no discussion is given as to whether they hold in this area of Biomechanics. The principle of conservation of mass is certainly valid, and the assumptions, of incompressibility of the fluid in the pores and of equilibrium of the porous medium as a whole, eminently plausible; other hypotheses, however, may be questionable or - at the very least - need to be assessed in a biological context. Can the brain, indeed, be treated as an isotropic and homogeneous medium as required by the theory? Does the assumed reversibility of the stress-strain relations under final equilibrium conditions hold? Does the linearity of the stress-strain relations (Hooke's law) apply?

Most biological materials do not obey this law exactly, for when they are suddenly strained, and this strain is held constant afterwards, the corresponding stress induced in them decreases with time - a phenomenon known as stress relaxation. It is therefore of interest to assess how important the viscoelastic properties of the brain are, and hence whether the theory of consolidation is applicable in the first place. Moreover, even if the validity of this theory is granted as a reasonable first approximation, the way it has been used by Nagashima **et al.** [17] is questionable. In particular, serious doubts must be raised against their use of the value $\nu = 0.4999$ for the Poisson ratio for white and gray matter. This value is practically equal to the theoretical maximum value that this parameter can take ($\nu = \frac{1}{2}$), and it is well known to be an upper limit never reached in real materials [10]. It implies that the material is perfectly incompressible, and has therefore the drastic consequence that the incremental volume strain $\epsilon = 0$. This, in conjunction with the fact that the CSF is also practically incompressible ($\beta' \simeq 0$) reduces the set of differential equations (38)-(40) to the Laplace equation $\nabla^2 \sigma = 0$, in which no variation with time is left. In view of this, the results presented by Nagashima **et al.** [17] in their Figure 4 - where the effective stress distribution at various times is displayed - are very difficult to understand.

Lastly, and most importantly, the question concerning the appropriate set of boundary conditions for the problem at hand must be squarely faced. As shown by many examples of application to soil mechanics and groundwater flow (see, for instance, Verruijt [24] and references therein) it may be possible to choose a geometry in which the set of Equations 38-40 along with the appropriate boundary conditions are amenable to approximate analytical solutions. Success in this direction would represent, in our opinion, an extremely important first step towards a full solution of the problem. Even if the model were extremely idealized, such analytical results would be valuable, for they would offer us formulas which could be used as a benchmark for testing the accuracy of a fully numerical solution. At the same time such theoretical studies would help to elucidate the role played by the physical parameters entering the model, such as the permeability and the porosity, which must ultimately be obtained from experimental measurements. It is not unreasonable to expect that many of these questions will be answered in the near future, for the problems facing researchers in this area of Biomechanics are not only important for their obvious clinical relevance, but also intellectually challenging and intriguing.

## APPENDICES

### A. Solution of Equation 8.

Before deriving the general solution of the basic DE describing the linear phenomenological models it is convenient to rewrite it is a simpler form. To this end we divide (**??**) through by $C_0$, and introduce the relaxation frequency $\nu_a = \dfrac{1}{C_0 R_a}$ and the equilibrium pressure $P_{\text{eq}}$ given by

(**??**) of the main text. Then (**??**) becomes simply

$$\frac{dP}{dt} + \nu_a P = \nu_a p_{\text{eq}} + \frac{1}{C_0} S(t). \tag{A1}$$

Next we solve this DE subject to the initial condition $P(0) = P_0$ using Laplace transforms. We define the transforms of $P(t)$ and $S(t)$ by the equations

$$\pi(s) = \int_0^\infty e^{-st} P(t) dt, \tag{A2}$$

$$\sigma(s) = \int_0^\infty e^{-st} S(t) dt, \tag{A3}$$

where $s$ is a real parameter having the physical dimensions of a frequency. Then multiplying (**??**) by $e^{-st}$, integrating with respect to time from zero to infinity, and using the above definitions and the known properties of the Laplace transform we obtain

$$s\pi(s) + \nu_a \pi(s) = P_0 + \nu_a P_{\text{eq}} \cdot \frac{1}{s} + \frac{1}{C_0} \sigma(s), \tag{A4}$$

which is no longer a DE, but just an algebraic equation for the transformed pressure $\pi$ whose solution is simply given by

$$\pi(s) = \frac{P_0}{s + \nu_a} + \frac{\nu_a P_{\text{eq}}}{s(s + \nu_a)} + \frac{1}{C_0} \frac{\sigma(s)}{s + \nu_a} \tag{A5}$$

Once the transformed function is explicitly known, as in (**??**), we can recover the original function $P(t)$ by application of the Inverse Laplace Transformation. For the four types of source functions listed in Table 1 the inverse Laplace transformation of (**??**) gives the solutions listed in the fifth column of the same Table, as the reader can easily verify.

**B. Solution of the Riccati Equation 21.**

It was mentioned in the main text of this paper that (**??**) is a well known DE, whose solution can be obtained by simple analytical means. All that is necessary is to perform the change of variables $P(t) = \frac{1}{p(t)}$ which reduces the equation to a linear differential equation (**??**). Since

$$\frac{dP}{dt} = -\frac{1}{p^2(t)} \frac{dp}{dt}, \tag{B1}$$

substitution into the original equation and multiplication of each resulting term by $-p^2(t)$ gives the result

$$\frac{dp}{dt} + \alpha(P_{\text{eq}} + R_a S(t))p = \alpha. \tag{B2}$$

The general solution of this first-order, linear DE is found by standard methods [6] and is given by

$$p(t) = e^{-\alpha P_{\text{eq}} t - \alpha R_a \int S(t)dt}\left(c + \alpha \int e^{\alpha P_{\text{eq}} t + \alpha R_a \int S(t)dt} dt\right), \tag{B3}$$

where $c$ is an arbitrary integration constant. Taking the reciprocal of this expression gives the desired solution for $P(t)$. This is now done for a few cases of interest:

**Case 1. No external source.** Mathematically, this means that $S(t) = 0$ for all $t$. Then (??) gives

$$p(t) = ce^{-\alpha P_{\text{eq}} t} + \frac{1}{P_{\text{eq}}} \tag{B4}$$

and the constant $c$ is determined from the initial condition $P(0) = P_0$, which implies that $p(0) = \frac{1}{P_0}$. Hence $c = \frac{1}{P_0} - \frac{1}{P_{\text{eq}}}$, and (??) becomes

$$p(t) = \frac{1}{P_{\text{eq}}} + \left(\frac{1}{P_0} - \frac{1}{P_{\text{eq}}}\right)e^{-\alpha P_{\text{eq}} t}. \tag{B5}$$

consequently the CSF pressure $P(t)$ is given by

$$P(t) = \frac{1}{\frac{1}{P_{\text{eq}}} + \left(\frac{1}{P_0} - \frac{1}{P_{\text{eq}}}\right)e^{-\alpha P_{\text{eq}} t}} = \frac{P_{\text{eq}}}{1 + \left(\frac{P_{\text{eq}}}{P_0} - 1\right)e^{-\alpha P_{\text{eq}} t}} \tag{B6}$$

**Case 2. Bolus injection (or withdrawal).**

When $S(t)$ is given by (??) we have

$$S(t) = \begin{cases} S_b, & \text{for } t_b \leq t \leq t_b + \Delta t_b \\ 0, & \text{otherwise} \end{cases}$$

and hence we have from (??) the following results:

(a) For $0 \leq t \leq t_b$ we have $S(t) = 0$, and so (??) reduces to (??). Imposing next the initial condition produces once again the solution (??), with the restriction, of course, that the time variable does not go beyond the value $t_b$.

(b) For $t_b \leq t \leq t_b + \Delta t_b$ we have $S(t) = S_b$, and so $\int S(t)dt = S_b t$. Consequently, (??) becomes

$$p(t) = ce^{-\alpha(P_{\text{eq}} + R_a S_b)t} + \frac{1}{P_{\text{eq}} + R_a S_b} \tag{B7}$$

where the arbitrary constant $c$ must be determined in such a way that $p(t)$ is continuous at $t = t_b$. When that is done, the result is

$$p(t) = \frac{1}{P_{\text{eq}}} \left[ \left(1 - \frac{P_{\text{eq}}}{P_b}\right) e^{\nu_b t_b} - \left(1 - \frac{P_{\text{eq}}}{P_0}\right) e^{\nu'_b t_b} \right] e^{-\nu_b t} + \frac{1}{P_b}, \tag{B8}$$

where, for simplification purposes, we have set

$$P_b = P_{\text{eq}} + R_a S_b, \tag{B9}$$

$$\nu_b = \alpha P_b, \tag{B10}$$

$$\nu'_b = \alpha R_a S_b. \tag{B11}$$

(c) For $t_b + \Delta t_b \leq t < \infty$ we have again $S(t) = 0$, and so the solution for $p(t)$ is again given by (??), where now the arbitrary constant $c$ must be determined so that $p$ is continuous at $t = t_b + \Delta t_b$. This gives

$$\begin{aligned} c = &\left(\frac{1}{P_b} - \frac{1}{P_{\text{eq}}}\right) e^{(\nu_b - \nu'_b)(t_b + \Delta t_b)} \\ &+ \frac{1}{P_{\text{eq}}} \left[ \left(1 - \frac{P_{\text{eq}}}{P_b}\right) e^{\nu_b t_b} - \left(1 - \frac{P_{\text{eq}}}{P_0}\right) e^{\nu'_b t_b} \right] e^{-\nu'_b (t_b + \Delta t_b)}, \end{aligned} \tag{B12}$$

which must be substituted into (??) and then the reciprocal taken in order to obtain $P(t)$; in other words

$$P(t) = \frac{P_{\text{eq}}}{c P_{\text{eq}} e^{(-\nu_b - \nu'_b)t} + 1}, \tag{B13}$$

where $c$ is given by (??).

ACKNOWLEDGMENTS

The research of S. Sivaloganathan and G. Tenti is supported in part by grants from the National Science and Engineering Research Council of Canada, and that of James M. Drake by the Appugliesi Hydrocephalus Research Fund.

**Figure 1.** Typical behavior of the CSF pressure as a function of time according to Eq. (11). The initial jump in the pressure to the value $P_0$ is assumed to be given.

**Figure 2.** A positive pressure pulse (bolus injection) at $t = t_0$ produces a sudden increase in the CSF pressure which relaxes to the equilibrium value thereafter. (Linear model.)

**Figure 3.** A sudden withdrawal of CSF at $t = t_0$ produces a sudden decrease in the CSF pressure which returns thereafter to the final equilibrium value. (Linear model.)

**Figure 4.** Typical behavior of the CSF pressure according to the linear model (Eq. 18) for the case of continuous infusion started at time $t_i$.

**Figure 5.** Comparison of the relaxation of the CSF pressure to the final equilibrium value as predicted by the linear (Eq. 11) and the non-linear (Eq. 32) models in the case of no external source.

**Figure 6.** Graph of the source function for the bolus injection case in the nonlinear model.

**Figure 7.** Recovery of the CSF pressure after an instantaneous bolus withdrawal according to the nonlinear phenomenological model (Eq. 36).

**Table 2.** Typical values of the parameters used in the numerical simulations.

| Symbol | Value | Units |
|---|---|---|
| $C_0$ | 0.004 | $cc/mmH_2O$ |
| $I_f^{(e)}$ | 0.078 | $cc/min$ |
| $k$ | 1 | $cc^{-1}$ |
| $P_0$ | 244 | $mmH_2O$ |
| $P_d$ | 70 | $mmH_2O$ |
| $\overline{P}$ | 40 | $mmH_2O$ |
| $P'$ | 150 | $mmH_2O$ |
| $P_{op}$ | 45 | $mmH_2O$ |
| $R_a$ | 600 | $mmH_2O/ml/min$ |
| $R_i$ | 600 | $mmH_2O/ml/min$ |
| $R_s$ | 60 | $mmH_2O/ml/min$ |
| $S_b$ | 40 | $cc/min$ |
| $S_i$ | 0.216 | $cc/min$ |
| $t_b, \Delta t_b, t_s$ | 5 | $min$ |
| $t_0, t_i$ | 4 | $min$ |

**Table 3.** Phenomenological Models of the Brain Compliance

| Model # | Functional Form of $C(P)$ | Corresponding Form of the Pressure-Volume Relation | Reference |
|---|---|---|---|
| 1 | $C(P) = C_0$ (constant) | $\Delta P = \frac{1}{C_0} \Delta V$ | Guinane (11) |
| 2 | $C(P) = \frac{1}{kP}$ | $P = P' e^{k\Delta V}$ | Marmarou **et al**. (15) |
| 3 | $C(P) = \frac{1}{k_1 P + k_2}$ ($k_1, k_2$) constant) | $P = \left(P' + \frac{k_2}{k_1}\right) e^{k_1 \Delta V} - \frac{k_2}{k_1}$ | Sklar **et al**. (19,20) |
| 4 | $C(P) = ae^{-bP}$ ($a, b$ positive constants) | $P = \ln\left(e^{-bP'} - \frac{b}{a}\Delta V\right)^{-1/b}$ | Lim **et al**. (14) |